\documentclass[12 pt]{article}

\usepackage{CJK}
\usepackage{amssymb}
\usepackage{bbm}
\usepackage{amsfonts}
\usepackage{mathrsfs}
\usepackage{graphicx}
\usepackage{amsmath}
\usepackage{amsthm}
\usepackage{xypic}
\usepackage{graphicx}
\usepackage{times}
\usepackage{geometry}
\usepackage{color}
\usepackage{indentfirst}
\usepackage{cases}
\usepackage{hyperref}
\usepackage{graphicx}
\usepackage{float}
\usepackage{amsthm}
\usepackage{amssymb}
\usepackage{amsmath}
\usepackage{mathrsfs}
\usepackage{indentfirst}
\usepackage{geometry}
\usepackage{authblk}
\usepackage{hyperref}

\hypersetup{
    colorlinks=true,
    linkcolor=blue,
    citecolor=blue
}
\newtheorem{Theorem}{\scshape  Theorem}[section]
\newtheorem{Lemma}[Theorem]{\scshape  Lemma}
\newtheorem{Corollary}[Theorem]{\scshape Corollary}

\newtheorem{Prop}[Theorem]{\scshape Proposition}

\newtheorem{Prob}[Theorem]{\scshape Problem}
\title{On oriented $m$-semiregular representations of finite groups about valency three}
 \author{
Songnian Xu\textsuperscript{a,}\thanks{Corresponding author. E-mail address: xsn131819@163.com},
\ Dein Wong\textsuperscript{a,}\thanks{Corresponding author. E-mail address: wongdein@163.com. Supported by the National Natural Science Foundation of China (No.12371025)},
\ Chi Zhang\textsuperscript{a,}\thanks{Supported by NSFC of China (No.12001526)},
\ Jinxing Zhao\textsuperscript{b,}\thanks{Supported by the National Natural Science Foundation of China (No.12161062) and Natural Science Foundation of Inner Mongolia (2025MS01022)}
\\
\textsuperscript{a}School of Mathematics, China University of Mining and Technology, Xuzhou, China.
\\
\textsuperscript{b}School of Mathematics Sciences, Inner Mongolia University, Hohhot, China.
}
\date{}

\begin{document}
\baselineskip 17pt

\title{On oriented $m$-semiregular representations of finite groups about valency three}

 \author{
Songnian Xu\textsuperscript{a,}\thanks{Corresponding author. E-mail address: xsn131819@163.com},
\ Dein Wong\textsuperscript{a,}\thanks{Corresponding author. E-mail address: wongdein@163.com. Supported by the National Natural Science Foundation of China (No.12371025)},
\ Chi Zhang\textsuperscript{a,}\thanks{Supported by NSFC of China (No.12001526)},
\ Jinxing Zhao\textsuperscript{b,}\thanks{Supported by the National Natural Science Foundation of China (No.12161062) and Natural Science Foundation of Inner Mongolia (2025MS01022)}
\\
\textsuperscript{a}School of Mathematics, China University of Mining and Technology, Xuzhou, China.
\\
\textsuperscript{b}School of Mathematics Sciences, Inner Mongolia University, Hohhot, China.
}
\date{}
\maketitle

\begin{abstract}
Let $G$ be a group and $m$ a positive integer. We say an $m$-Cayley digraph $\Gamma$ over $G$ is a digraph that admits a group of automorphisms isomorphic to $G$ acting semiregularly on the vertex set with $m$ orbits.
The digraph $\Sigma$ is $k$-regular if there exists a non-negative integer $k$ such that every vertex has out-valency and in-valency equal to $k$. All digraphs considered in this paper are regular.
We say that $G$ admits an oriented $m$-semiregular representation (abbreviated as OmSR) if there exists a regular $m$-Cayley digraph $\Gamma$ over $G$ such that $\Gamma$ is oriented and its automorphism group is isomorphic to $G$.
In particular, an O1SR is called an ORR. Xia et al. \cite{x2} provided a classification of finite simple groups admitting an ORR of valency 2. Furthermore, in 2022, Du et al. \cite{du2} proved that most finite simple groups admit an OmSR of valency 2 for $m \geq 2$, except for a few exceptional cases. In this paper, we classify the finite groups generated by at most two elements that admit an OmSR of valency 3 for $m \geq 2$.
\end{abstract}

\let\thefootnoteorig\thefootnote
\renewcommand{\thefootnote}{\empty}

\section{Introduction}
A \textit{directed graph} (or \textit{digraph}) $\Gamma$ is defined as an ordered pair $(V(\Gamma), A(\Gamma))$, where $V(\Gamma)$ denotes a non-empty vertex set and $ A(\Gamma)$ is a subset of $V(\Gamma) \times V(\Gamma)$. The elements of $V$ and $A$ are referred to as \textit{vertices} and \textit{arcs}, respectively. For any arc $(u, v) \in A$, $v$ is called an \textit{out-neighbor} of $u$, while $u$ is termed an \textit{in-neighbor} of $v$. 
The $out$-$valency$ (resp. $in$-$valency$) of a vertex $v \in V(\Gamma)$ is the number of its out-neighbors (resp. in-neighbors). 
A digraph $\Gamma$ is said to be \textit{$k$-regular} if there exists a non-negative integer $k$ such that each vertex in $\Gamma$ has both out-degree and in-degree equal to $k$. All digraphs discussed in this paper are assumed to be regular.
Given a subset $X \subseteq V(\Gamma)$, the \textit{subdigraph} of $\Gamma$ induced by $X$ is defined as $\Gamma[X] := (X, A(\Gamma) \cap (X \times X))$, often abbreviated as $[X]$. A digraph $\Gamma$ is a \textit{graph} if its arc set $A(\Gamma)$ is symmetric, meaning that $(u, v) \in A(\Gamma)$ implies $(v, u) \in A(\Gamma)$. Conversely, a digraph is \textit{oriented} if for every pair of distinct vertices $u$ and $v$, at most one of the ordered pairs $(u, v)$ or $(v, u)$ is an arc in $A(\Gamma)$.

Let $\Gamma$ be a digraph and $\omega \in V(\Gamma)$. 
An \textit{automorphism} of $\Gamma$ is a permutation $\sigma$ of $V(\Gamma)$ that leaves the arc set $A(\Gamma)$ invariant, i.e.,
\[
(x^\sigma, y^\sigma) \in A(\Gamma) \iff (x, y) \in A(\Gamma).
\]
The set of all automorphisms of $\Gamma$, equipped with the composition operation, forms the \textit{full automorphism group} of $\Gamma$, denoted by $\operatorname{Aut}(\Gamma)$. For a subgroup $G \leq \operatorname{Aut}(\Gamma)$, let $G_\omega$ denote the \textit{stabilizer} of $\omega$ in $G$, i.e., the subgroup of $G$ that fixes $\omega$. We say that $G$ acts \textit{semiregularly} on $V(\Gamma)$ if $G_\omega = 1$ for every $\omega \in V(\Gamma)$, and \textit{regularly} if it acts semiregularly and transitively.

Let $G$ be a finite group and $S$ be a subset of $G$. The $Cayley$ $digraph$ $\mathrm{Cay}(G,S)$ has vertex set $G$ and arc set $\{(g,sg) \mid g \in G, s \in S\}$. The graph is a \emph{Cayley graph} (resp. an oriented Cayley digraph) if and only if the connection set $S$ satisfies the symmetric condition $S = S^{-1}$ (resp. $S \cap S^{-1} = \emptyset$).
The concept of Cayley digraphs can be nicely generalized to $m$-$Cayley$ $digraphs$ where regular actions are replaced with semiregular actions. An $m$-Cayley (di)graph $\Gamma$ over a finite group $G$ is defined as a (di)graph which has a semiregular group of automorphisms isomorphic to $G$ with $m$ orbits on its vertex set.

We say a finite group $G$ admits a \emph{graphical regular representation} (GRR), \emph{digraphical regular representation} (DRR), or \emph{oriented regular representation} (ORR) if there exists a Cayley graph $\Gamma = \mathrm{Cay}(G, S)$, Cayley digraph $\Gamma = \mathrm{Cay}(G, S)$, or oriented Cayley digraph $\Gamma = \mathrm{Cay}(G, S)$ over $G$, respectively, such that $\mathrm{Aut}(\Gamma) \cong G$.
The GRR and DRR problem asks which groups admit such representations. Babai \cite{bab1} proved that except for $Q_8$, $\mathbb{Z}_2^2$, $\mathbb{Z}_2^3$, $\mathbb{Z}_2^4$, and $\mathbb{Z}_3^2$, every group admits a DRR. Spiga \cite{mor,mor1,spi1} classified groups admitting ORRs. While GRR implies DRR, the converse fails. The GRR turned out to be much more difficult to handle and, after a long series of partial results by various authors \cite{he,im,im1,im2,no,no1,wa}, the classification was completed by Godsil in \cite{god1}.

We say a finite group $G$ admits:
\begin{itemize}
\item a \emph{graphical m-semiregular representation} (GmSR) if there exists a regular $m$-Cayley graph $\Gamma$ over $G$ with $\mathrm{Aut}(\Gamma) \cong G$;
\item a \emph{digraphical m-semiregular representation} (DmSR) if there exists a regular $m$-Cayley digraph $\Gamma$ over $G$ with $\mathrm{Aut}(\Gamma) \cong G$;
\item an \emph{oriented m-semiregular representation} (OmSR) if there exists a regular oriented $m$-Cayley digraph $\Gamma$ over $G$ with $\mathrm{Aut}(\Gamma) \cong G$.
\end{itemize}
Note that G1SR, D1SR and O1SR coincide with GRR, DRR and ORR respectively. Groups admitting GmSR, DmSR or OmSR for all $m$ were classified in \cite{du1,du1'}.

\textbf{OmSR with prescribed valency $k$}: Unlike unrestricted GmSR, DmSR, and OmSR, the problem of classifying groups that admit GmSR, DmSR, or OmSR with a specified valency remains largely unsolved.
For digraphs, the minimum possible value is 2 when \( k \) is even, and 1 when \( k \) is odd. That said, a connected digraph of valency 1 is merely a directed cycle when \( k = 1 \). As such, the smallest non-trivial case is valency 2 for even \( k \), and valency 3 for odd \( k \).
As for the valency 2 scenario, Verret and Xia \cite{x2} classified finite simple groups with an ORR (i.e., O1SR) of valency 2, showing that every simple group of order at least 5 has such an ORR.
Du et al. \cite{du1'} provided a full classification of OmSR for all finite groups, and their later work \cite{du2} described OmSR of valency 2 specifically for finite groups generated by at most two elements.
These findings naturally prompt an investigation into classifying OmSR of valency 3 for finite simple groups.
By the classification theorem of finite simple groups, all such groups can be generated by at most two elements.
Thus, in this paper, we examine OmSR of valency 3 for groups generated by at most two elements.

Of course, there are many more articles on the automorphism groups of graphs, which we will not mention one by one here; interested readers may refer to \cite{du3,du4,du5,ou1,ou2,ou3,ou4,ou5} and so on.
Below, we present the main results of this paper.

\begin{Theorem}
Let $G = \langle x \rangle\neq 1$ be a finite cyclic group and $m \geq 2$ be a positive integer. Then $G$ admits an OmSR of valency 3, except when either:
\begin{itemize}
\item  $o(x)=2$ and $m\leq 4$ ,
\item  $o(x) = 3$ or $4$ and $m = 2$.
\end{itemize}
\end{Theorem}

\begin{Theorem}
Let $G=\langle x,y\rangle\neq\langle x\rangle$ be a finite group and $m\geq 2$ an integer. Then $G$ admits an OmSR of valency 3, except when $G\cong\mathbb{Z}_2^2$ and $m=2$.
\end{Theorem}

The following corollary is an immediate consequence of Theorems~1.1 and 1.2.

\begin{Corollary}

Let $G$ be a non-trivial finite simple group $m\geq 2$ an integer. Then $G$ admits an OmSR of valency 3, except for when:
\begin{enumerate}
    \item[(i)] $G \cong \mathbb{Z}_2$ and $m\leq4$,
    \item[(ii)]  $G \cong \mathbb{Z}_3$ and $m=2$.
\end{enumerate}

\end{Corollary}

In the case where $G = 1$ (the trivial group), we have yet to resolve the classification of OmSRs of valency 3 for $G$. This problem is actually equivalent to the following basic question: For each integer $m$, does there exist a regular digraph $\Gamma$ of order $m$ and valency 3 with a trivial automorphism group?

To end this section, we propose the following problems:

\begin{Prob}
 For every integer $m$, does there exist a regular digraph $\Gamma$ of order $m$ and valency 3 with trivial automorphism group?
\end{Prob}

\section{Preliminaries and notations}

Let $m$ be a positive integer and $G$ a group. To simplify notation in this paper, we write the element $(g,i)$ of the Cartesian product $G \times \{0, \ldots, m-1\}$ as $g_i$. We often identify $\{0, \ldots, m-1\}$ with $\mathbb{Z}_m$, the integers modulo $m$.
Recall that an $m$-Cayley digraph of a finite group $G$ is a digraph admitting a semiregular group of automorphisms isomorphic to $G$ with exactly $m$ orbits on its vertex set. 
We now present a more specific definition of the $m$-Cayley digraph.

For each $i \in \mathbb{Z}_m$, define $G_i := \{g_i \mid g \in G\}$. Analogous to classical Cayley digraphs, an $m$-Cayley digraph can be viewed as the digraph
\[
\Gamma = \mathrm{Cay}(G, T_{i,j} : i, j \in \mathbb{Z}_m)
\]
with:
\begin{itemize}
\item Vertex set: $G \times \mathbb{Z}_m = \bigcup_{i \in \mathbb{Z}_m} G_i$,
\item Arc set: $\bigcup_{i,j \in \mathbb{Z}_m} \{(g_i, (tg)_j) \mid t \in T_{i,j}\}$,
\end{itemize}
where $T_{i,j} \subseteq G$ for all $i,j \in \mathbb{Z}_m$.

When $T_{i,j}\cap T^{-1}_{j,i} = \emptyset$ for all $i \in \mathbb{Z}_m$, we call $\Gamma$ an \emph{oriented $m$-Cayley digraph}.
For any $g \in G$, the right multiplication map $R(g)$, defined by $R(g): x_i \mapsto (xg)_i$ for all $x_i \in G_i$ and $i \in \mathbb{Z}_m$, is an automorphism of $\Gamma$.
$R(G)=\{R(g) \mid g \in G\}$ is isomorphic to $G$ and it is a semiregular group of automorphisms of $\Gamma$ with  $G_i$ as  orbits.

For a digraph $\Gamma$ and any vertex $x \in V(\Gamma)$, we define:
\begin{itemize}
   \item $ \Gamma^{+0}(x) = \{x\}$,
  \item $\Gamma^{+1}(x) = \Gamma^+(x)$ is the set of out-neighbors of $x$ in $\Gamma$ \text{ and } $\Gamma^+[x]=\Gamma^+(x)\cup \{x\}$,
    \item $\Gamma^{+2}(x)$ represents the union of out-neighbors of all vertices in $\Gamma^+(x)$,
    \item  $\Gamma^{+k}(x) = \bigcup_{y \in \Gamma^{+(k-1)}(x)} \Gamma^+(y) \quad \text{for any integer } k \geq 1$.
 \end{itemize}

Finally, we end this section with Proposition 2.1, which is a crucial tool for our proof that a finite group admits an OmSR.

\begin{Prop}\cite[3.1]{du3}
Let $m$ be a positive integer at least $2$ and let $G$ be a finite group.
For any $i, j\in \mathbb{Z}_m$, let $T_{i,j}\subseteq G$ and let $\Gamma=Cay(G,T_{i,j}: i,j\in \mathbb{Z}_m)$ be a connected $m$-cayley digraph over $G$.
For $A=Aut(\Gamma)$, if $A$ fixes $G_i$ setwise for all $i\in \mathbb{Z}_m$ and there exist $u_0\in G_0,u_1\in G_1,\ldots,u_{m-1}\in G_{m-1}$ such that $A_{u_i}$ fixes $\Gamma^{+}(u_i)$ pointwise for all $i\in \mathbb{Z}_m$, then $A=R(G)$.
\end{Prop}

\section{Proof of Theorems}
\begin{Lemma}
Let $G=\langle x\rangle\cong \mathbb{Z}_2$ be a finite group and $m\geq 2$ be a positive integer. 
Then $G$ admits an OmSR of valency 3 if and only if $m\geq 5$.
\end{Lemma}

\begin{proof}
Let $\Gamma=\mathrm{Cay}(G_i,T_{i,j}: i,j\in \mathbb{Z}_m)$ be a connected oriented $m$-Cayley digraph and $A=\mathrm{Aut}(\Gamma)$.

Since $\Gamma$ is an oriented digraph with $o(x)=2$, it immediately follows that $T_{i,i}=\emptyset$. Moreover, when $T_{i,j}=G$, the condition $T_{i,j}\cap T^{-1}_{j,i}=\emptyset$ for $i,j\in \mathbb{Z}_m$ implies $T_{j,i}=\emptyset$.

\noindent\textbf{Case 1: $m=2$}\\
Since $|G|=2$ and $T_{i,i}=\emptyset$ for all $i\in\mathbb{Z}_m$, the resulting digraph $\Gamma$ cannot be a 3-regular digraph.

\noindent\textbf{Case 2: $m=3$}\\
As $\Gamma$ is 3-regular and oriented, either $T_{0,1}=G$ or $T_{0,2}=G$ by $T_{i,i}=\emptyset$ and $o(x)=2$ for $i\in \mathbb{Z}_3$. 
Without loss of generality, assume $T_{0,1}=G$. 
Then $T_{1,2}=G$ since $T_{1,0}=T_{1,1}=\emptyset$. 
Similarly, we have $T_{2,0}=G$.
Finally, we obtain $T_{0,1} = T_{1,2} = T_{2,3} = G$, with all other $T_{i,j} = \emptyset$ for $i \in \mathbb{Z}_3$, and the resulting digraph $\Gamma$ cannot satisfy the 3-regularity condition.

\noindent\textbf{Case 3: $m=4$}\\
If there exists no $T_{i,j} = G$ for $0 \leq i \neq j \leq 3$, then since $\Gamma$ is 3-regular and $T_{i,j} \cap T^{-1}_{j,i} = \emptyset$, we must have $T_{i,j} \cup T_{j,i} = G$. 
However, using MAGMA or Mathematica, we find that $A>R(G)$. 
To keep the paper concise, we provide the method for computing the automorphism group of a given graph Gamma using Mathematica or MAGMA in the final section.

If there exists some $T_{i,j}=G$, assume $T_{0,1}=G$. Then from $T_{1,0}=\emptyset$ and $\Gamma$ being 3-regular, we obtain either $T_{1,2}=G$ or $T_{1,3}=G$.

First consider $T_{1,2}=G$. Then $T_{2,1}=\emptyset$, and by 3-regularity, $|T_{1,3}|=1$. This leads to two subcases: $T_{2,0}=G$ or $T_{2,3}=G$.
\begin{itemize}
\item If $T_{2,0}=G$, then  $T_{0,2}=\emptyset, |T_{2,3}|=1$ and $|T_{0,3}|=1$. Consequently, $|T_{3,0}|=|T_{3,1}|=|T_{3,2}|=1$. 
At this point, we have determined all $T_{i,j}$ for $i,j \in \mathbb{Z}_4$:
\begin{align*}
T_{0,1} &= T_{0,2} = T_{2,0} = G, \\
|T_{i,3}| &= |T_{3,i}| = 1, T_{i,3} \cup T_{3,i} = G \quad \text{for } i \in \mathbb{Z}_4\setminus\{3\}, \\
T_{i,j} &= \emptyset \quad \text{for all other cases}.
\end{align*}
Using Mathematica or MAGMA, we find that for the graph $\Gamma$ determined by these $T_{i,j}$, $A \neq R(G)$.

\item If $T_{2,3}=G$, then $T_{3,2}=\emptyset$ and by 3-regularity, $|T_{2,0}|=1$, $T_{3,0}=G$, and $|T_{3,1}|=1$. Similarly, we get $|T_{0,2}|=1$. 
We now obtain all $T_{i,j}$ for $i,j \in \mathbb{Z}_4$:
\begin{align*}
T_{0,1} &= T_{1,2} = T_{2,3} = T_{3,0} = G, \\
|T_{0,2}| &= |T_{2,0}| = |T_{1,3}| = |T_{3,1}| = 1 \text{ and } T_{0,2} \cup T_{2,0} = T_{1,3} \cup T_{3,1} = G, \\
T_{i,j} &= \varnothing \quad \text{for all other cases.}\\
\text{Again, computation shows $A\neq R(G)$.}
\end{align*}
\end{itemize}
The case $T_{1,3}=G$ similarly leads to $A\neq R(G)$.

\noindent\textbf{Case 4: $5\leq m\leq 10$}\\
Define $T_{i,j}=\emptyset$ except for the following:
\begin{itemize}
\item For $m=5$: \\
$T_{0,1}=T_{0,3}=T_{2,0}=T_{2,1}=T_{2,4}=T_{3,1}=T_{3,2}=T_{3,4}=T_{4,1}=\{1\}$, \\ $T_{1,0}=T_{1,2}=T_{1,3}=T_{1,4}=T_{4,3}=T_{4,2}=\{x\}$.
\item For $m=6$: \\
$T_{1,0}=T_{0,3}=T_{0,5}=T_{1,4}=T_{2,0}=T_{2,3}=T_{2,5}=T_{3,1}=T_{4,3}=T_{5,4}=1$, \\ $T_{1,0}=T_{1,2}=T_{3,2}=T_{3,4}=T_{4,1}=T_{4,5}=T_{5,0}=T_{5,2}=\{x\}$.
\item For $m=7$: \\
$T_{0,1}=T_{2,1}=T_{2,0}=T_{2,5}=T_{3,1}=T_{3,2}=T_{3,4}=T_{4,5}=T_{5,6}=T_{6,4}=T_{6,0}=\{1\}$, \\ $T_{0,6}=T_{0,3}=T_{1,0}=T_{1,2}=T_{1,3}=T_{4,3}=T_{4,6}=T_{5,2}=T_{5,4}=T_{6,5}=\{x\}$.
\item For $m=8$: \\
$T_{0,1}=T_{1,5}=T_{2,0}=T_{2,3}=T_{2,5}=T_{3,1}=T_{3,4}=T_{4,6}=T_{5,6}=T_{6,7}=T_{7,4}=T_{7,0}=\{1\}$,\\ $T_{0,3}=T_{0,7}=T_{1,0}=T_{1,2}=T_{3,2}=T_{4,3}=T_{4,7}=T_{5,1}=T_{5,2}=T_{6,4}=T_{6,5}=T_{7,6}=\{x\}$.
\item For $m=9$: \\ $T_{0,1}=T_{2,0}=T_{2,1}=T_{2,5}=T_{3,1}=T_{3,2}=T_{3,4}=T_{4,6}=T_{5,6}=T_{5,7}=T_{6,7}=T_{7,8}=T_{8,0}=T_{8,1}=\{1\}$, \\ $T_{0,3}=T_{0,8}=T_{1,0}=T_{1,2}=T_{1,3}=T_{4,3}=T_{4,8}=T_{5,2}=T_{6,4}=T_{6,5}=T_{7,5}=T_{7,6}=T_{8,7}=\{x\}$.
\item For $m=10$: \\ $T_{0,1}=T_{1,4}=T_{2,0}=T_{2,3}=T_{2,5}=T_{3,1}=T_{3,4}=T_{5,6}=T_{5,7}=T_{6,7}=T_{6,8}=T_{7,8}=T_{8,9}=T_{9,4}=T_{9,0}=\{1\}$,\\
 $T_{0,3}=T_{0,9}=T_{1,0}=T_{1,2}=T_{3,2}=T_{4,1}=T_{4,3}=T_{4,9}=T_{5,2}=T_{6,5}=T_{7,5}=T_{7,6}=T_{8,6}=T_{8,7}=T_{9,8}=\{x\}$.
\end{itemize}

Clearly, for $5 \leq m \leq 10$, the graph $\Gamma$ obtained from the $T_{i,j}$ defined above is a 3-regular oriented graph.
Using Mathematica or MAGMA, we verify that $A=R(G)\cong G$ for $5\leq m\leq 10$.

\noindent\textbf{Case 5: $m\geq 11$}\\
Define:
\begin{itemize}
\item $T_{i,i+1}=\{1\}$ for $i\in \mathbb{Z}_m$
\item $T_{i,i-2}=\{1,x\}$ except $T_{2,5}=T_{7,0}=\{1,x\}$ for $i\in \mathbb{Z}_m\backslash\{2,7\}$
\item $T_{i,j}=\emptyset$ otherwise
\end{itemize}
The digraph $\Gamma$ is shown in Fig.1 and is clearly 3-regular oriented graph. We now prove $A=R(G)\cong G$.
\begin{figure}[H]
  \centering
  \includegraphics[width=1.0\linewidth]{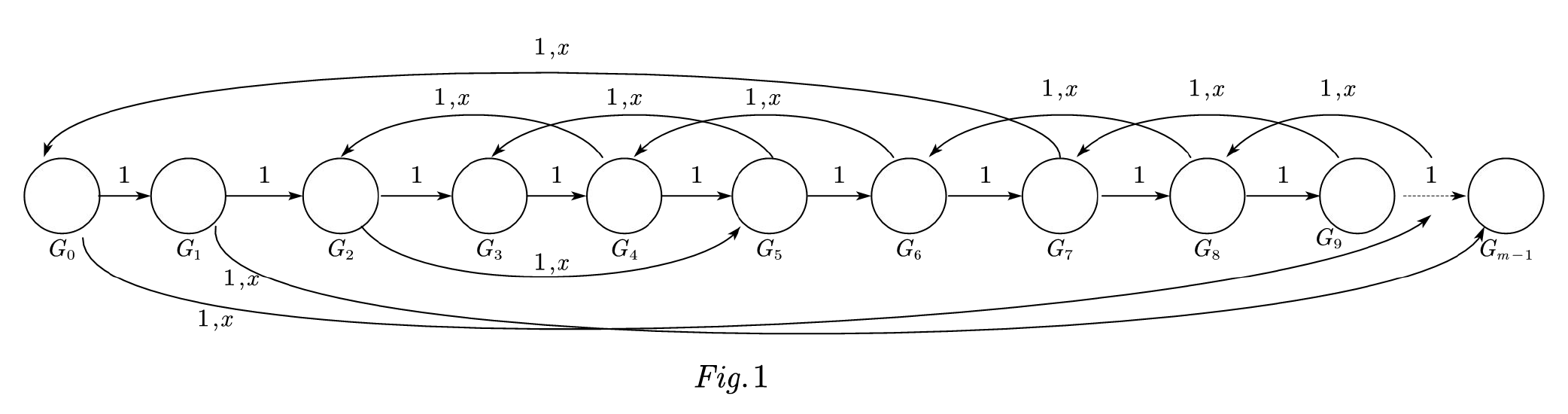}
\end{figure}

We first prove that $A$ fixes each $G_i$ setwise for $i\in \mathbb{Z}_m$.

Observe that $(g_i,g_{i+1},g_{i+2})$ forms a oriented cycle, except when $i=0$ or $5$ (where $i\in\mathbb{Z}_m$ and $g\in G$). Therefore, $A$ fixes both $H_1=G_0\cup G_1\cup G_2\cup G_5\cup G_6\cup G_7$ and $H_2=G_3\cup G_4\cup (\bigcup_{i\geq 8}G_i)$ setwise.

Note that $[H_2]=[G_3\cup G_4]+[\bigcup_{i\geq 8}G_i]$ (for two graphs $\Gamma_1$ and $\Gamma_2$, we denote by $\Gamma_1 + \Gamma_2$ their \emph{disjoint union}), where:
\begin{itemize}
\item $[G_3\cup G_4]$ is disconnected because $T_{3,4}=\{1\}$.
\item $[\bigcup_{i\geq 8}G_i]$ is connected via $T_{10,8}=\{1,x\}$.
\end{itemize}
Thus $A$ fixes $G_3\cup G_4$ setwise.

Moreover:
\begin{itemize}
\item Vertices in $G_3$ have out-valency 2 in $H_1$ (via $T_{3,1}=\{1,x\}$).
\item Vertices in $G_4$ have out-valency 3 in $H_1$ (via $T_{4,5}=\{1\}$ and $T_{4,2}=\{1,x\}$).
\end{itemize}
Therefore $A$ fixes both $G_3$ and $G_4$ setwise individually.

Further analysis shows:
\begin{itemize}
\item Since $T_{3,i}=\emptyset$ for $i\neq 1,4$ and $A$ fixes both $G_3$ and $G_4$ setwise, it follows that $A$ fixes $G_1$ setwise.
\item The oriented 3-cycle $(g_1,g_2,g_3)$ is unique (no $(g_1,g_i,g_3)$ exists for $i\neq 2$). Since $A$ fixes both $G_1$ and $G_3$ setwise, it must necessarily fix $G_2$ setwise.
\item Since $T_{4,i}=\emptyset$ for $i\neq 2,5$ and $A$ fixes both $G_2$ and $G_4$ setwise, it follows that $A$ must fix $G_5$ setwise.
\item Since $T_{5,i}=\emptyset$ for $i\neq 3,6$ and $A$ fixes both $G_3$ and $G_5$ setwise, it follows that $A$ must fix $G_6$ setwise.
\item From the conditions $T_{6,4}=\{1,x\}$, $T_{6,7}=\{1\}$, and $T_{6,i}=\emptyset$ for $i\neq 4,7$, and since $A$ fixes both $G_4$ and $G_6$ setwise, we conclude that $A$ must fix $G_7$ setwise.

\end{itemize}

Combining with $A$ fixing $H_1=G_0\cup\cdots\cup G_7$, we conclude $A$ fixes $G_0$ setwise. Thus $A$ fixes $G_i$ setwise for all $0\leq i\leq 7$.

For $i\geq 8$:
\begin{itemize}
\item From $T_{7,8}={1}$, $T_{7,0}={1,x}$, and $T_{7,i}=\emptyset$ for all other $i$, since $A$ fixes both $G_0$ and $G_7$ setwise, it follows that $A$ also fixes $G_8$ setwise.
\item Similarly, $A$ fixes $G_i$ for $i\geq 9$ via $T_{i,i-2}=\{1,x\}$, $T_{i,i+1}=1$, and $T_{i,j}=\emptyset$ otherwise.
\end{itemize}

Therefore, $A$ fixes $G_i$ setwise for all $i\in\mathbb{Z}_m$.
By the Frattini argument, we have $A = R(G)A_{g_i}$, which means:
\begin{align}
|A_{g_i}| = |A_{h_j}| = |A/R(G)|,  \tag{1}
\end{align}
for all $g,h \in G$ and $i,j \in \mathbb{Z}_m$.
We now prove $A_{1_0}=1$. Since $T_{i,i+1}=1$ for $i\in\mathbb{Z}_m$, $A_{1_i}$ fixes $1_{i+1}$, so $A_{1_i}=A_{1_{i+1}}$  by (1). Similarly, $A_{x_i}=A_{x_{i+1}}$ for $i\in\mathbb{Z}_m$. Thus:
\begin{align}
A_{1_0}=A_{1_1}=\cdots=A_{1_{m-1}}, \tag{2}\\
A_{x_0}=A_{x_1}=\cdots=A_{x_{m-1}}. \tag{3}
\end{align}

Since $\Gamma^+(1_0)=\{1_1,1_{m-2},x_{m-2}\}$, $A_{1_0}$ fixes $x_{m-2}$ by (2), so $A_{1_0}=A_{x_{m-2}}$ by (1). Combining with (2) and (3):
\begin{equation}
A_{1_0}=A_{1_1}=\cdots=A_{1_{m-1}}=A_{x_0}=A_{x_1}=\cdots=A_{x_{m-1}}. \tag{4}
\end{equation}

Now, we prove $A_{x_0}=A_{x^2_0}$ using the fact that $A_{1_0}=A_{x_0}$; consider $\sigma\in A_{x_0}$. Since:
\[ 1_0^{\rho_x\sigma\rho_{x^{-1}}} = x_0^{\sigma\rho_{x^{-1}}} = x_0^{\rho_{x^{-1}}} = 1_0, \]
we have $\rho_x\sigma\rho_{x^{-1}}\in A_{1_0}=A_{x_0}$. Thus:
\[ x_0^{\rho_x\sigma\rho_{x^{-1}}} = x_0 \Rightarrow (x_0^2)^\sigma = x_0^2, \]
so $\sigma\in A_{x_0^2}$. Therefore $A_{x_0}=A_{x_0^2}$ by equal orders.

From (4) and the conjugate action of $R(G)$, we obtain:
\begin{equation}
A_{x_0} = A_{x^2_0}= \cdots= A_{1_0} =A_{x_1}= \cdots = A_{1_1}=\cdots=A_{x_{m-1}} = A_{x^2_{m-1}}= \cdots= A_{1_{m-1}}. \tag{5}
\end{equation}

Since $G$ is cyclic and $o(x)=2$, (4) already implies $A_{1_0}=1$, i.e., $A=R(G)\cong G$. For subsequent arguments, we derive (5) from (4), though this step will be omitted in later discussions.
\end{proof}

\begin{Lemma}
Let $G = \langle x \rangle$ be a finite cyclic group and $m \geq 2$ be a positive integer with $o(x) \geq 3$. Then $G$ admits an OmSR of valency 3, except when  $m = 2$ and $o(x) = 3$ or $o(x) = 4$.

\end{Lemma}

\begin{proof}
Let $\Gamma=\mathrm{Cay}(G_i,T_{i,j}: i,j\in \mathbb{Z}_m)$ be a connected oriented $m$-Cayley digraph and $A=\mathrm{Aut}(\Gamma)$.

\subsection*{Case 1: $m=2$}

\subsubsection*{Case 1.1: $o(x)=3$}
Since $G=\{1,x,x^2\}$ and $\Gamma$ is an oriented digraph, we have $1\notin T_{i,i}$ for $i\in \mathbb{Z}_2$.
If $T_{0,0}=\{x,x^2\}$, then $T_{0,0}\cap T^{-1}_{0,0}=T_{0,0}\neq \emptyset$, which leads to a contradiction. Therefore $|T_{0,0}|\leq 1$. Similarly, $|T_{1,1}|\leq 1$.
Since $\Gamma$ is required to be 3-regular, we must have $|T_{1,0}|,|T_{0,1}|\geq 2$. In this case, we obtain $T_{0,1}\cap T^{-1}_{1,0}\neq \emptyset$ because $o(x)=3$.
In conclusion, when $o(x)=3$, the group $G=\langle x\rangle$ does not admit an O2SR of valency 3.

\subsubsection*{Case 1.2: $o(x)=4$}
Here $G=\{1,x,x^2,x^3\}$. Clearly, we also have $1,x^2=x^{-2}\notin T_{i,i}$ for $i\in \mathbb{Z}_2$.
If $|T_{0,0}|\geq 2$, then by the condition $o(x)=4$, we must have $T_{0,0}=\{x,x^3\}$ and consequently $T_{0,0}\cap T^{-1}_{0,0}=T_{0,0}$, which leads to a contradiction.
Thus $|T_{0,0}|\leq 1$, and similarly $|T_{1,1}|\leq 1$.
If $T_{0,0}=\emptyset$, then $|T_{0,1}|=3$ since $\Gamma$ is a 3-regular digraph. Given $|T_{1,1}|$, we have $|T_{1,0}|\geq 2$, which necessarily implies $T_{1,0}\cap T^{-1}_{0,1}\neq \emptyset$ because $o(x)=4$, leading to a contradiction. Therefore $|T_{0,0}|=1$, and similarly $|T_{1,1}|=1$.
Combining with $x^2,1\notin T_{i,i}$ for $i\in \mathbb{Z}_2$, we have $T_{i,i}\subseteq \{\{x\},\{x^3\}\}$.

Since $T_{0,1}\cap T^{-1}_{1,0}=\emptyset$ and $|T_{0,1}|=|T_{1,0}|=2$ (because $\Gamma$ is an oriented 3-regular digraph), it follows that $T_{1,0},T_{0,1}\in \{\{1,x\},\{x,x^2\}\}$ with $T_{0,1}\neq T_{1,0}$.

At this point, we have determined all possible choices for $T_{i,j}$ where $i,j\in \mathbb{Z}_2$. However, using Mathematica or MAGMA, we find that for all $\Gamma$ satisfying these conditions, $A>R(G)\cong G$. Therefore, when $o(x)=4$, the group $G=\langle x\rangle$ does not admit an O2SR of valency 3.

\subsubsection*{Case 1.3: $o(x)\geq 5$}
Let $T_{0,0}=\{x,x^2\}$, $T_{0,1}=\{1\}$; $T_{1,1}=\{x^{-1},x^{-2}\}$, $T_{1,0}=\{x\}$. Then $\Gamma$ is clearly an oriented 3-regular digraph. We now prove that $A=R(G)\cong G$.
We have:
\begin{align*}
\Gamma^{+}(1_0) &= \{x_0, x^2_0, 1_1\},
\Gamma^{+}(1_1) = \{x^{-1}_1, x^{-2}_1, x_0\}; \\
\Gamma^{+2}(1_0) &= \Gamma^{+}(x_0) \cup \Gamma^{+}(x^2_0) \cup \Gamma^{+}(1_1) \\
&= \{x^2_0, x^3_0, x_1\} \cup \{x^3_0, x^4_0, x^2_1\} \cup \{x^{-1}_1, x^{-2}_1, x_0\} \\
&= \{x_0, x^2_0, x^3_0, x^4_0, x_1, x^{-1}_1, x^{2}_1, x^{-2}_1\}, \\
\Gamma^{+2}(1_1) &= \Gamma^{+}(x^{-1}_1) \cup \Gamma^{+}(x^{-2}_1) \cup \Gamma^{+}(x_0) \\
&= \{x^{-2}_1, x^{-3}_1, 1_0\} \cup \{x^{-3}_1, x^{-4}_1, x^{-1}_0\} \cup \{x^{2}_0, x^{3}_0, x_1\} \\
&= \{1_0, x^{-1}_0, x^2_0, x^3_0, x_1, x^{-2}_1, x^{-3}_1, x^{-4}_1\}.
\end{align*}
\begin{figure}[H]
  \centering
  \includegraphics[width=0.4\linewidth]{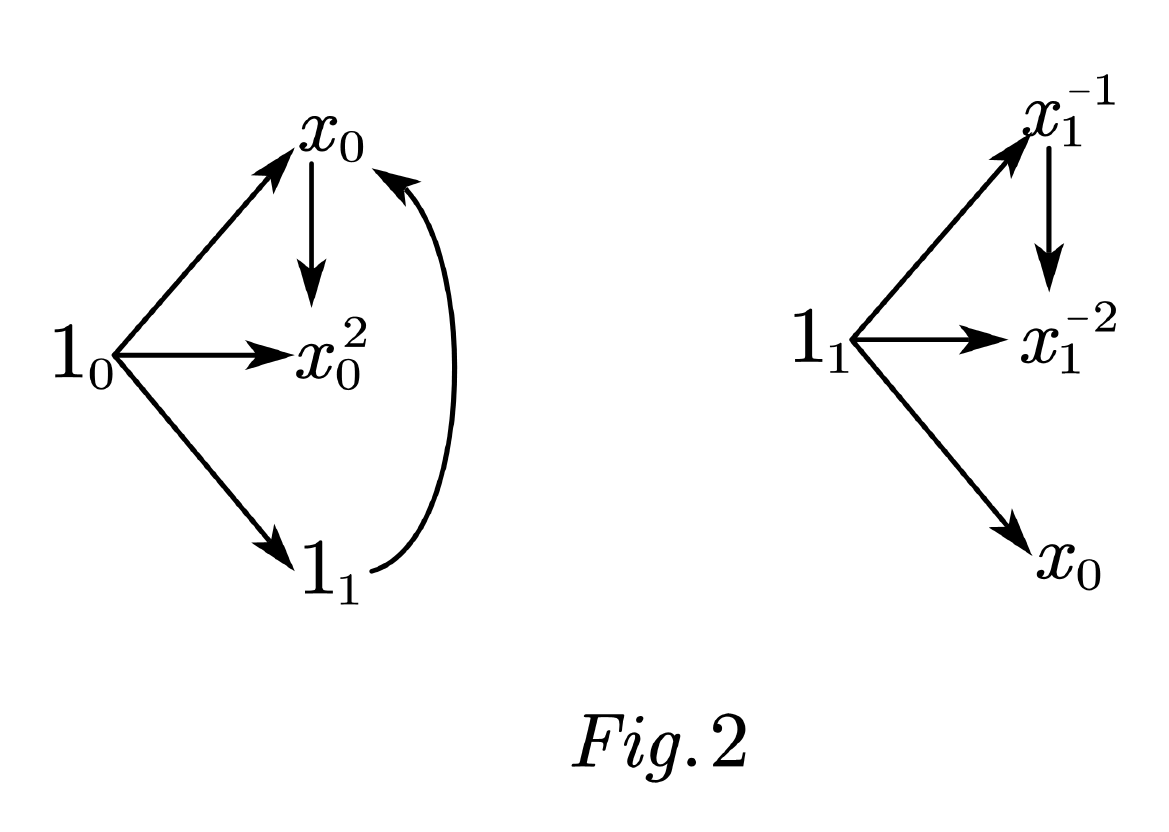}
\end{figure}

The graph $[\Gamma^{+}(1_i)]$ is shown in Fig.2 for each $i \in \mathbb{Z}_2$.
Note that in $[\Gamma^{+}(1_0)]$ there are two arcs $(x_0,x^2_0)$ and $(1_0,x_0)$, while in $[\Gamma^{+}(1_1)]$ there is only one arc $(x^{-1}_1,x^{-2}_1)$. Therefore, $A$ fixes both $G_0$ and $G_1$ setwise.

By the Frattini argument, we have:
\[ A = A_{1_0}R(G) = A_{x_0}R(G) = A_{1_1}R(G) = A_{x_1}R(G). \]

Now we prove that $A_{1_0} = A_{x_0} = A_{1_1} = A_{x_1}$:
\begin{itemize}
\item Since $T_{0,1}=\{1\}$, we have $A_{1_0}=A_{1_1}$ and $A_{x_0}=A_{x_1}$.
\item Since $T_{1,0}=\{x\}$, we have $A_{1_0}=A_{x_1}$.
\end{itemize}
In conclusion, we obtain $A_{1_0}=A_{x_0}=A_{1_1}=A_{x_1}$.

Furthermore, we have:
\[ A_{x_0} = A_{x^2_0} = \cdots = A_{1_0} = A_{x_1} = A_{x^2_1} = \cdots = A_{1_1}, \]
and since $G$ is cyclic, it follows that $A_{1_0}=1$, i.e., $A=R(G)\cong G$.

\subsection*{Case 2: $m\geq 3$}
Define:
\begin{align*}
T_{0,0} &= \{x\},
T_{i,i} = \{x^{-1}\} \text{ for } i\in \mathbb{Z}_m\backslash\{0\}, \\
T_{m-1,1}&=\{x,x^{-1}\}, T_{i,i+1} = \{1,x^{-1}\} \text{ for } i\in \mathbb{Z}_m\backslash\{m-1\}, \\
T_{i,j} &= \emptyset \text{ otherwise.}
\end{align*}

\begin{figure}[H]
  \centering
  \includegraphics[width=1.1\linewidth]{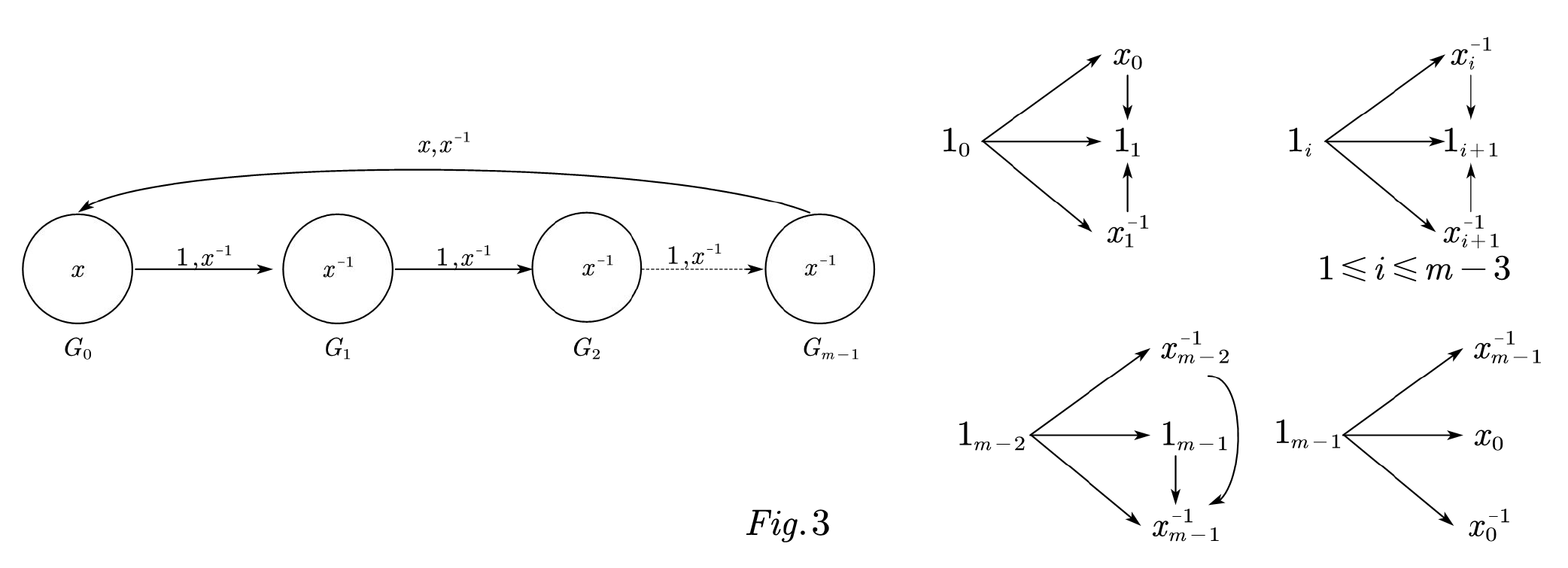}
\end{figure}
Then the graph $\Gamma$ is shown in Fig.3 and $\Gamma$ has the following properties:\\
$\Gamma^{+}(1_0) = \{x_0, 1_1, x^{-1}_1\}, \\
\Gamma^{+}(1_i) = \{x^{-1}_i, 1_{i+1}, x^{-1}_{i+1}\} \text{ for } 1\leq i\leq m-2, \\
\Gamma^{+}(1_{m-1}) = \{x^{-1}_{m-1}, x_0, x^{-1}_0\}$;
\\
\\
$\Gamma^{+2}(1_0) = \Gamma^{+}(x_0) \cup \Gamma^{+}(1_1) \cup \Gamma^{+}(x^{-1}_1) \\
= \{x^2_0, x_1, 1_1\} \cup \{x^{-1}_1, 1_2, x^{-1}_2\} \cup \{x^{-2}_1, x^{-1}_2, x^{-2}_2\} \\
= \{x^2_0, 1_1, x_1, x^{-1}_1, x^{-2}_1, 1_2, x^{-1}_2, x^{-2}_2\}$;
\\
\\
$\Gamma^{+2}(1_i) = \Gamma^{+}(x^{-1}_i) \cup \Gamma^{+}(1_{i+1}) \cup \Gamma^{+}(x^{-1}_{i+1}) \\
= \{x^{-2}_i, x^{-1}_{i+1}, x^{-2}_{i+1}\} \cup \{x^{-1}_{i+1}, 1_{i+2}, x^{-1}_{i+2}\}  \cup \{x^{-2}_{i+1}, x^{-1}_{i+2}, x^{-2}_{i+2}\} \\
= \{x^{-2}_i, x^{-1}_{i+1}, x^{-2}_{i+1}, 1_{i+2}, x^{-1}_{i+2}, x^{-2}_{i+2}\}
 \text{ for } 1\leq i\leq m-3$;
 \\
 \\
$\Gamma^{+2}(1_{m-2}) = \Gamma^{+}(x^{-1}_{m-2}) \cup \Gamma^{+}(1_{m-1}) \cup \Gamma^{+}(x^{-1}_{m-1}) \\
= \{x^{-2}_{m-2}, x^{-1}_{m-1}, x^{-2}_{m-1}\} \cup \{x^{-1}_{m-1}, x_{0}, x^{-1}_{0}\}  \cup \{x^{-2}_{m-1}, 1_{0}, x^{-2}_{0}\} \\
= \{1_0, x_0, x^{-1}_0, x^{-2}_0, x^{-2}_{m-2}, x^{-1}_{m-1}, x^{-2}_{m-1}\}$;
\\
\\
$\Gamma^{+2}(1_{m-1}) = \Gamma^{+}(x^{-1}_{m-1}) \cup \Gamma^{+}(x_{0}) \cup \Gamma^{+}(x^{-1}_{0}) \\
= \{x^{-2}_{m-1}, 1_{0}, x^{-2}_{0}\} \cup \{x^{2}_{0}, x_{1}, 1_{1}\} \cup \{1_{0}, x^{-1}_{1}, x^{-2}_{1}\} \\
= \begin{cases}
\{1_0, x^{-2}_0, x^2_0, 1_1, x_1, x^{-1}_1, x^{-2}_1, x^{-2}_{m-1}\}  \text{ for } o(x)\geq 5 \\
\{1_0, x^2_0, 1_1, x_1, x^{-1}_1, x^{-2}_1, x^{-2}_{m-1}\}  \text{ for } o(x)=4 \\
\{1_0, x_0, x^2_0, 1_1, x_1, x^{-1}_1, x^{-2}_{m-1}\}  \text{ for } o(x)=3
\end{cases}
$

For clarity, we illustrate $[\Gamma^{+}[1_i]]$ for $i\in \mathbb{Z}_m$ as shown in Fig.3.

We observe that:
\begin{itemize}
\item In $[\Gamma^{+}(1_0)]$, there are two arcs $(x_0,1_1)$ and $(1_1,x^{-1}_1)$.
\item In $[\Gamma^{+}(1_i)]$ for $1\leq i\leq m-2$, there are two arcs $(x^{-1}_i,x^{-1}_{i+1})$ and $(1_{i+1},x^{-1}_{i+1})$.
\item In $[\Gamma^{+}(1_{m-1})]$, when $o(x)=3$ there is one arc $(x^{-1}_{m-1},x^{-2}_0)=(x^{-1}_{m-1},x_0)$, while for $o(x)>3$, $[\Gamma^{+}(1_{m-1})]$ forms an isolated set.
\end{itemize}

Therefore, $A$ fixes $G_{m-1}$ setwise. Moreover, since $T_{i,j}=\emptyset$ for $j\neq i$ or $i+1$, it follows that $A$ fixes $G_i$ for all $i\in \mathbb{Z}_m$.

By the Frattini argument, we have $A=A_{g_i}R(G)$ for $g\in G$ and $i\in \mathbb{Z}_m$, since $A$ fixes each $G_i$ setwise.

Consequently:
\begin{itemize}
\item $A_{1_0}$ fixes both $\{x_0\}$ and $\{1_1,x^{-1}_1\}$ setwise, i.e., $A_{1_0}=A_{x_0}$;
\item $A_{1_i}$ fixes both $\{x^{-1}_i\}$ and $\{1_{i+1},x^{-1}_{i+1}\}$ setwise for $1\leq i\leq m-2$, i.e., $A_{1_i}=A_{x^{-1}_i}$;
\item $A_{1_{m-1}}$ fixes both $\{x^{-1}_{m-1}\}$ and $\{x_0,x^{-1}_0\}$ setwise, i.e., $A_{1_{m-1}}=A_{x^{-1}_{m-1}}$.
\end{itemize}

Furthermore, we note that:
\begin{itemize}
\item In $[\Gamma^{+}(1_0)]$, $(x_0,1_1)$ is an arc while $(x_0,x^{-1}_1)$ is not, so $A_{1_0}$ fixes $1_1$ (since $A_{1_0}$ fix $x_0$), i.e., $A_{1_0}=A_{1_1}$.
\item In $[\Gamma^{+}(1_i)]$, $(x^{-1}_i,x^{-1}_{i+1})$ is an arc while $(x^{-1}_i,1_{i+1})$ is not, so $A_{1_i}$ fixes $1_{i+1}$ (since $A_{1_i}$ fixes $x^{-1}_i$) for $1\leq i\leq m-2$, hence $A_{1_1}=A_{1_2}=\cdots=A_{1_{m-1}}$.
\end{itemize}

In conclusion, we have:
\[ A_{1_0} = A_{x_0} = A_{1_1} = A_{x^{-1}_1} = \cdots = A_{1_{m-1}} = A_{x^{-1}_{m-1}}. \]

By the conjugate action of $R(G)$ on $G_i$, we establish that:
\[ A_{x_0} = A_{x^2_0} = \cdots = A_{1_0} = A_{x_1} = A_{x^2_1} = \cdots = A_{1_1} = \cdots = A_{x_{m-1}} = A_{x^2_{m-1}} = \cdots = A_{1_{m-1}}. \]

Since $G$ is cyclic, we have $A_{1_0}=1$, i.e., $A=R(G)\cong G$.
\end{proof}

Theorem 1.1 follows directly from Lemmas 3.1 and 3.2.

We now consider the case where $G=\langle x,y\rangle\neq\langle x\rangle$, dividing our analysis into two subcases: $G\cong\mathbb{Z}_2^2$ and $G\ncong\mathbb{Z}_2^2$.

\begin{Lemma}
Let $G=\langle x,y\rangle\cong \mathbb{Z}_2^2$ and $m\geq 2$ an integer. Then $G$ admits an OmSR of valency 3, except when $m=2$.
\end{Lemma}

\begin{proof}

Let $\Gamma=\mathrm{Cay}(G_i,T_{i,j}: i,j\in \mathbb{Z}_m)$ be a connected oriented $m$-Cayley digraph and $A=\mathrm{Aut}(\Gamma)$.
\subsubsection*{Case 1: $m=2$} For $m=2$, \cite{du1'} shows that $G$ does not admit an O2SR of valency 3.

\subsubsection*{Case 2: $m=3$ or $4$}

For $m=3$, define:
\begin{align*}
T_{0,1} &= \{1\}, \quad T_{0,2} = \{xy,x\}; \\
T_{1,0} &= \{x,y\}, \quad T_{1,2} = \{1\}; \\
T_{2,0} &= \{y\}, \quad T_{2,1} = \{x,y\}, \\
T_{i,j} &= \emptyset \text{ for other } i,j\in\mathbb{Z}_3.
\end{align*}

For $m=4$, define:
\begin{align*}
T_{0,1} &= \{1,x\}, \quad T_{0,2} = \{xy\}; \\
T_{1,0} &= \{y\}, \quad T_{1,3} = \{1,y\}; \\
T_{2,0} &= \{1,x\}, \quad T_{2,3} = \{x\}; \\
T_{3,1} &= \{x\}, \quad T_{3,2} = \{1,y\}, \\
T_{i,j} &= \emptyset \text{ for other } i,j\in\mathbb{Z}_4.
\end{align*}

These constructions yield 3-regular oriented digraphs $\Gamma$. For $m=3$ or $4$, the digraph $\Gamma$ constructed from the $T_{i,j}$ defined above is clearly a 3-regular oriented graph. 
Using Mathematica or MAGMA, we verify that $A=R(G)\cong G$. 
For clarity, we provide a detailed proof for $m=3$.

\noindent\textbf{Detailed Proof for $m=3$:}

Neighborhoods:
\begin{align*}
\Gamma^{+}(1_0) &= \{1_1, x_2, (xy)_2\}, \\
\Gamma^{+}(1_1) &= \{1_2, x_0, y_0\}, \\
\Gamma^{+}(1_2) &= \{y_0, x_1, y_1\}.
\end{align*}

Second neighborhoods:
\begin{align*}
\Gamma^{+2}(1_0) &= \Gamma^{+}(1_1)\cup \Gamma^{+}(x_2)\cup \Gamma^{+}((xy)_2)\\
&=\{1_2,x_0,y_0\}\cup \{(xy)_0,1_1,(xy)_1\}\cup \{x_0,y_1,x_1\}\\
&=\{1_2, x_0, y_0, (xy)_0, (xy)_1, y_1, x_1\}, \\
\Gamma^{+2}(1_1) &= \Gamma^{+}(1_2)\cup \Gamma^{+}(x_0)\cup \Gamma^{+}((y_0)\\
&=\{y_0,x_1,y_1\}\cup \{x_1,1_2,y_2\}\cup \{y_1,(xy)_2,x_2\}\\
&=\{y_0, x_1, y_1, 1_2, y_2, (xy)_2, x_2\}, \\
\Gamma^{+2}(1_2) &= \Gamma^{+}(y_0)\cup \Gamma^{+}(x_1)\cup \Gamma^{+}((y_1)\\
&=\{y_1,(xy)_2,x_2\}\cup \{x_2,1_0,(xy)_0\}\cup \{y_2,(xy)_0,1_0\}\\
&= \{y_1, (xy)_2, x_2, 1_0, (xy)_0, y_2\}.
\end{align*}

\begin{figure}[H]
  \centering
  \includegraphics[width=1.0\linewidth]{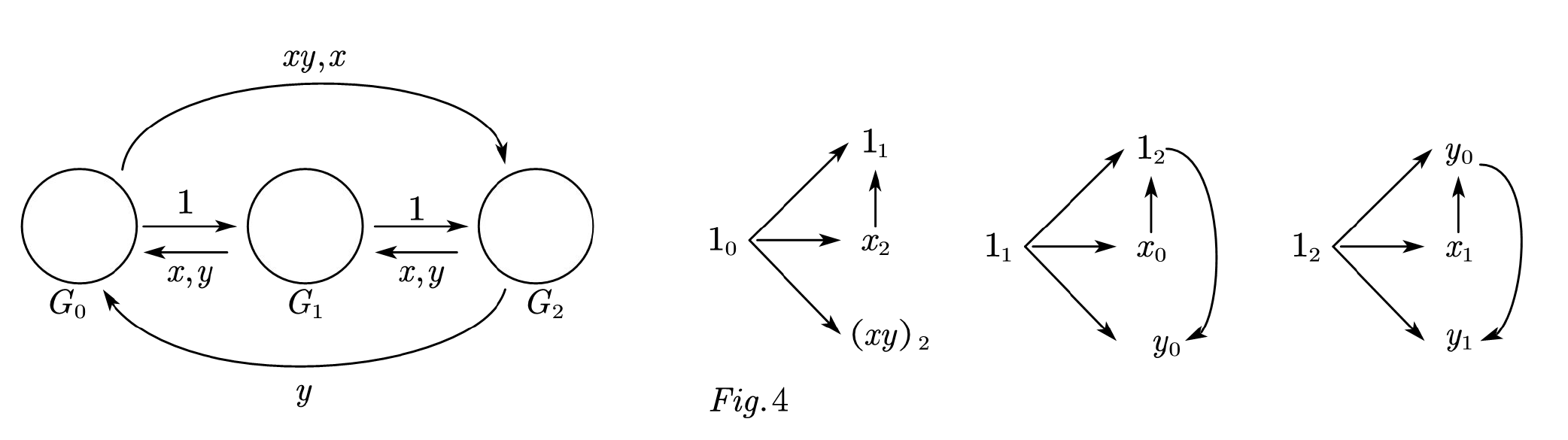}
\end{figure}

Fig. 4 shows $\Gamma$ and $[\Gamma^{+}[1_i]$ for $i\in \mathbb{Z}_3$. We can see:
\begin{itemize}
\item $[\Gamma^{+}(1_0)]$ and $[\Gamma^{+}(1_2)]$ contain isolated vertices;
\item $[\Gamma^{+}(1_1)]$ has no isolated vertices.
\end{itemize}

Thus $A$ fixes $G_1$ and $G_0\cup G_2$ setwise. Moreover:
\begin{itemize}
\item Vertices in $G_0$ have out-valency 1 in $G_1$ (via $T_{0,1}=\{1\}$),
\item Vertices in $G_2$ have out-valency 2 in $G_1$ (via $T_{2,1}=\{x,y\}$).
\end{itemize}
Therefore $A$ fixes both $G_0$ and $G_2$ setwise.
Now we prove $A_{1_i}$ fixes $\Gamma^{+}(1_i)$ pointwise for $i \in \mathbb{Z}_3$. 
Since $A$ fixes each $G_i$ setwise ($i \in \mathbb{Z}_3$), we have:
\begin{itemize}
\item $A_{1_0}$ fixes $\{1_1\}$ and $\{x_2,(xy)_2\}$ setwise,
\item $A_{1_1}$ fixes $\{1_2\}$ and $\{x_0,y_0\}$ setwise,
\item $A_{1_2}$ fixes $\{y_0\}$ and $\{x_1,y_1\}$ setwise.
\end{itemize}

Since in $[\Gamma^{+}(1_0)]$, $(x_2, 1_1)$ is an arc and $((xy)_2, 1_1)$ is not an arc, we have that $A_{1_0}$ fixes $x_2$ and $(xy)_2$ because $A_{1_0}$ fixes $1_1$. That is, $A_{1_0}$ fixes $\Gamma^{+}(1_0)$ pointwise. For $i = 1, 2$, using the same argument, we also have that $A_{1_i}$ fixes $\Gamma^{+}(1_i)$ pointwise. By Proposition 2.1, $A=R(G)\cong G$.

\subsubsection*{Case 3: $m\geq5$}
Define the sets $T_{ij}\subseteq G$:
\begin{align*}
T_{0,1} &= \{y\}, \\
T_{i,i+1} &= \{1\} \text{ for } i\in \mathbb{Z}_m\backslash\{0\}, \\
T_{i,i-2} &= \{1,x\} \text{ for } i\in \mathbb{Z}_m, \\
T_{i,j} &= \emptyset \text{ otherwise.}
\end{align*}

\begin{figure}[H]
  \centering
  \includegraphics[width=0.8\linewidth]{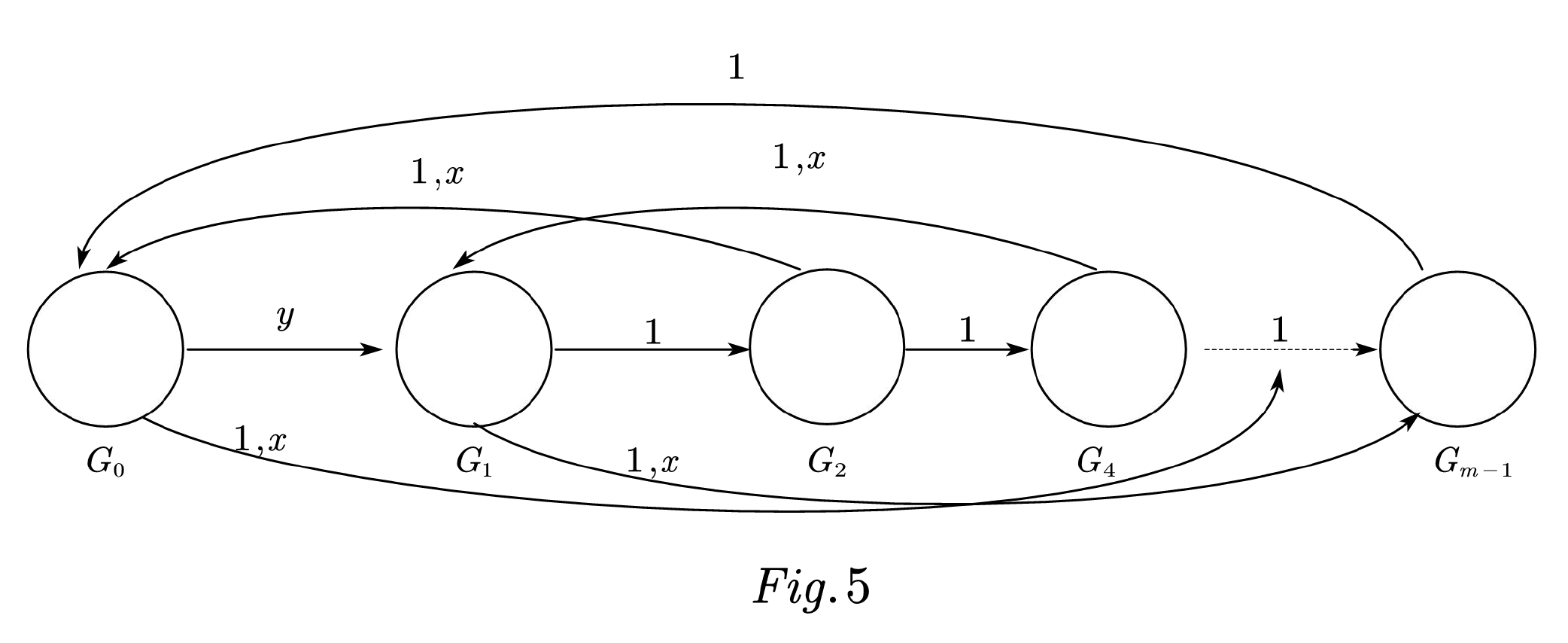}
\end{figure}

The resulting digraph $\Gamma$ (shown in Fig. 5) is clearly a 3-regular oriented digraph. We now prove $A=R(G)\cong G$.

From Fig. 5, we easily observe that $(g_i, g_{i+1}, g_{i+2})$ forms a oriented 3-cycle for all $i$ except $i = 0$ and $i = m-1$. Thus, $A$ fixes $H_1 = G_{m-1} \cup G_0 \cup G_1 \cup G_2$ and $H_2 = V(\Gamma) \backslash H_1$ setwise. We now prove that $A$ fixes $G_i$ for $i \in \mathbb{Z}_m$.

When $m = 5$, we observe that vertices in $G_1$ have in-valency 2 in $H_2 = G_3$, vertices in $G_0$ and $G_2$ have in-valency 0 in $H_2$, and vertices in $G_{m-1} = G_4$ have in-valency 1 in $H_2 $. Thus, $A$ fixes $G_1$ setwise.

When $m > 5$, we observe that vertices in $G_1$ and $G_2$ have in-valency 2 in $H_2$ by $T_{3,1} = T_{4,2} = \{1, x\}$, vertices in $G_0$ have in-valency 0 in $H_2$ by $T_{i,0} = \emptyset$ for $i \neq 2, m-1$, and vertices in $G_{m-1}$ have in-valency 1 in $H_2$ by $T_{m-2,m-1} = \{1\}$. Thus, $A$ fixes $G_1 \cup G_2$ setwise.
On the other hand, vertices in $G_1$ have out-valency 0 in $H_2$ by $T_{1,i}=\emptyset$ for $i\neq2,m-1$, while vertices in $G_2$ have out-valency 1 in $H_2$ by $T_{2,3}=\{1\}$. Therefore, $A$ fixes $G_1$ setwise.

In summary, for $m \geq 5$, we have that $A$ fixes $G_1$. We now prove that $A$ fixes $G_i$ for $i \in \mathbb{Z}_m$.

$T_{1,i} = \emptyset$ for $i \neq 2, m-1$ implies that $A$ fixes $G_2 \cup G_{m-1}$ setwise. Moreover, since $|T_{1,2}| = 1$ and $|T_{1,m-1}| = 2$, it follows that $A$ fixes both $G_2$ and $G_{m-1}$ setwise because $A$ fixes $G_1$ setwise. 

$T_{2,i} = \emptyset$ for $i \neq 0, 3$ implies that $A$ fixes $G_0 \cup G_3$ setwise. Further, combining with $|T_{2,3}| = 1$ and $|T_{2,0}| = 2$, we have that $A$ fixes both $G_0$ and $G_3$ setwise because $A$ fixes $G_2$ setwise.
By similar arguments, we can show that $A$ fixes $G_i$ setwise for all $i \in \mathbb{Z}_m$.
We now prove that $A_{1_i}$ fixes $\Gamma^{+}(1_i)$ for $i \in \mathbb{Z}_m$.

Since $T_{i,i+1}=\{1\}$ for $i\in\mathbb{Z}_m\backslash\{0\}$, we have:
\[ A_{1_1} = A_{1_2} = \cdots = A_{1_{m-1}} = A_{1_0}. \]

Since $T_{0,1}=\{y\}$, we have $A_{1_0}=A_{y_1}$, leading to:
\[ A_{1_0} = A_{1_1} = \cdots = A_{1_{m-1}} = A_{y_1} \quad \text{(6)}. \]

Neighborhood properties:
\begin{itemize}
\item $\Gamma^{+}(1_0) = \{y_1, 1_{m-2}, x_{m-2}\}$.
\item $\Gamma^{+}(1_i) = \{1_{i+1}, 1_{i-2}, x_{i-2}\}$ for $i\in\mathbb{Z}_m\backslash\{0\}$.
\end{itemize}

From equation (6), we conclude that $A_{1_i}$ fixes $\Gamma^{+}(1_i)$ pointwise for all $i\in\mathbb{Z}_m$. Therefore, by Proposition 2.1, $A=R(G)\cong G$.
\end{proof}

\begin{Lemma}
Let $G=\langle x,y\rangle\ncong \mathbb{Z}_2^2 $ be a finite group and $m\geq 2$ an integer. Then $G$ admits an OmSR of valency 3.
\end{Lemma}

\begin{proof}
Since $G\ncong \mathbb{Z}_2^2$, we have either $o(x)\geq 3$ or $o(y)\geq 3$. Without loss of generality, assume $o(x)\geq 3$.

Define the connection sets:
\begin{align*}
T_{i,i} &= \{x\} \text{ for } i\in \mathbb{Z}_m \\
T_{m-1,1} &= \{x,y\}, \\
T_{i,i+1} &= \{1,x\} \text{ for } i\in \mathbb{Z}_m\backslash\{m-1\}, \\
T_{i,j} &= \emptyset \text{ otherwise.}
\end{align*}

\begin{figure}[H]
  \centering
  \includegraphics[width=0.9\linewidth]{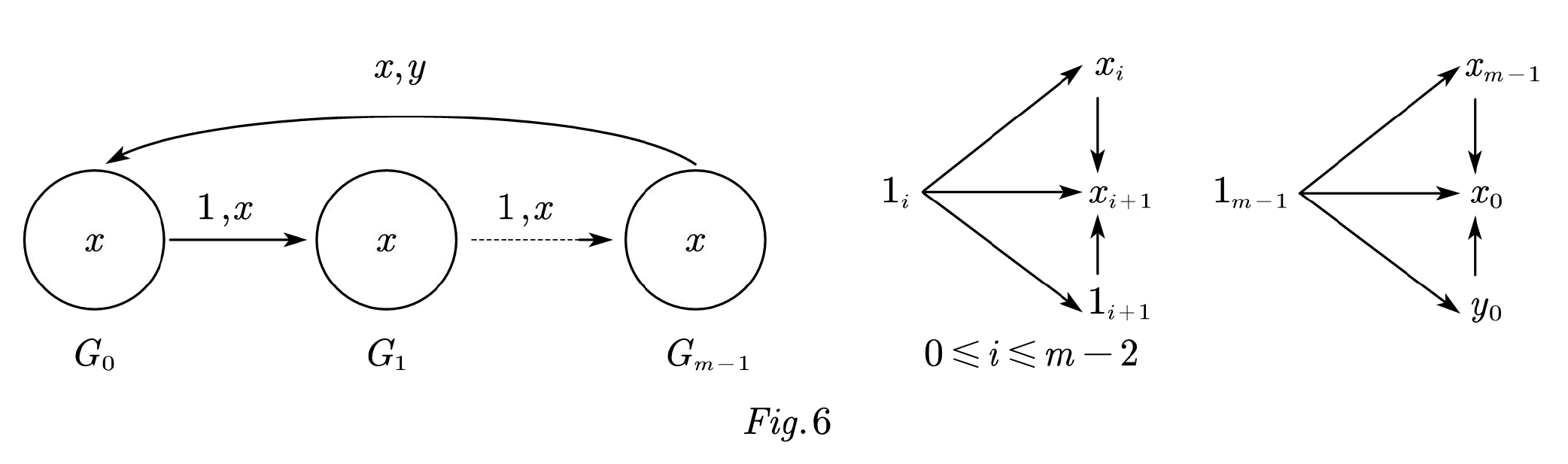}
\end{figure}
The resulting digraph $\Gamma$ (shown in Fig. 6) is clearly a 3-regular oriented digraph. We now prove $A=R(G)$.First, we show $A$ fixes each $G_i$ setwise for $i\in \mathbb{Z}_m$.
\begin{align*}
\Gamma^{+}(1_i) &= \{x_i, 1_{i+1}, x_{i+1}\} \text{ for } 0\leq i\leq m-2, \Gamma^{+}(1_{m-1}) = \{x_{m-1}, x_0, y_0\}.
\end{align*}

Second neighborhoods:
\begin{align*}
\Gamma^{+2}(1_i) &= \Gamma^{+}(x_i) \cup \Gamma^{+}(1_{i+1}) \cup \Gamma^{+}(x_{i+1}) \\
&= \{x^{2}_i, x_{i+1}, x^{2}_{i+1}\} \cup \{x_{i+1}, 1_{i+2}, x_{i+2}\} \cup \{x^{2}_{i+1}, x_{i+2}, x^{2}_{i+2}\} \\
&= \{x^2_i, x_{i+1}, x^2_{i+1}, 1_{i+2}, x_{i+2}, x^{2}_{i+2}\} \quad \text{for } 0 \leq i \leq m-3\\
\Gamma^{+2}(1_{m-2}) &= \Gamma^{+}(x_{m-2}) \cup \Gamma^{+}(1_{m-1}) \cup \Gamma^{+}(x_{m-1}) \\
&= \{x^{2}_{m-2}, x_{m-1}, x^{2}_{m-1}\} \cup \{x_{m-1}, x_{0}, y_{0}\} \cup \{x^{2}_{m-1}, x^2_{0}, (yx)_{0}\} \\
&= \{x^2_{m-2}, x_{m-1}, x^2_{m-1}, x_{0}, x^2_{0}, y_{0}, (yx)_0\}\\
\Gamma^{+2}(1_{m-1}) &= \Gamma^{+}(x_{m-1}) \cup \Gamma^{+}(x_{0}) \cup \Gamma^{+}(y_{0}) \\
&= \{x^{2}_{m-1}, x^2_{0}, (yx)_{0}\} \cup \{x^2_{0}, x_{1}, x^2_{1}\} \cup \{xy_{0}, y_1, (xy)_1\} \\
&= \{x^2_0, (yx)_0, (xy)_{0}, x_{1}, x^2_{1}, y_{1}, (xy)_1, x^2_{m-1}\}
\end{align*}

From Fig. 6, we observe that in $[\Gamma^{+}(1_i)]$, there are two arcs $(x_i, x_{i+1})$ and $(1_{i+1}, x_{i+1})$ for $0 \leq i \leq m-2$; in $[\Gamma^{+}(1_{m-1})]$, there exists an arc $(x_0, x_1)$ when $m = 2$, but no arcs exist when $m \geq 3$. Therefore, $A$ fixes $G_{m-1}$ setwise. Since $T_{m-1,i} = \emptyset$ for $i \neq m-1, 0$, $A$ fixes $G_0$ setwise. Similarly, $A$ fixes $G_i$ for all $i \in \mathbb{Z}_m$ because $T_{i,j} = \emptyset$ for $j \neq i, i+1$.

We show that $A_{1_i}$ fixes $\Gamma^{+}(1_i)$ pointwise for all $i \in \mathbb{Z}_m$:

For $0 \leq i \leq m-2$, $A_{1_i}$ fixes both $\{x_i\}$ and $\{1_{i+1}, x_{i+1}\}$ setwise because $A$ fixes $G_i$ setwise. On the other hand, in $[\Gamma^{+}(1_i)]$, $(x_i, x_{i+1})$ is an arc while $(x_i, 1_{i+1})$ is not. Therefore, $A_{1_i}$ fixes $\Gamma^{+}(1_i)$ pointwise for $0 \leq i \leq m-2$ because $A_{1_i}$ fixes $x_i$.
Now, we only need to prove that $A_{m-1}$ fixes $\Gamma^{+}(1_{m-1})$.
By the Frattini argument, we have:
\begin{equation}
A_{1_0} = A_{1_1} = \cdots = A_{1_{m-1}} = A_{x_0} = A_{x_1} = \cdots = A_{x_{m-1}}    \tag{7}. 
\end{equation}

Since $\Gamma^{+}(1_{m-1})=\{x_0,y_0,x_{m-1}\}$, equation (7) implies $A_{1_{m-1}}$ fixes $\Gamma^{+}(1_{m-1})$ pointwise.
Therefore, $A_{1_i}$ fixes $\Gamma^{+}(1_i)$ pointwise for all $i\in\mathbb{Z}_m$, and by Proposition 2.1, we conclude that $A=R(G)\cong G$.
\end{proof}
Theorem 1.2 follows directly from Lemmas 3.3 and 3.4.

\section{Computing Automorphism Groups in Mathematica}
The following is the Mathematica code for computing the automorphism group of a graph:

\begin{verbatim}
(* === Mathematica Code to Compute Graph Automorphism Group === *)

(* Step 1: Define your graph \(\Gamma\) *)
(* Example 1: Cycle graph with 4 vertices *)
\[Gamma] = CycleGraph[4]; 

(* Example 2: Custom graph (uncomment to use) *)
(* \[Gamma] = Graph[{1->2, 2->3, 3->1, 1->4, 4->3}, 
VertexLabels->"Name"] *)

(* Step 2: Compute automorphism group *)
autGroup = GraphAutomorphismGroup[\[Gamma]];

(* Step 3: Display results *)
Print["Automorphism Group Order: ", GroupOrder[autGroup]];
Print["Group Generators: ", GroupGenerators[autGroup]];
Print["All Group Elements: ", GroupElements[autGroup]];

(* Optional visualization *)
GraphPlot[\[Gamma], VertexLabeling -> True]

\end{verbatim}

To facilitate readers' understanding and application of this program, we demonstrate how to compute the automorphism group of a graph using Mathematica, based on Case 1.2 of Lemma 3.2 in this paper.

In Case 1.2 of Lemma 3.2, we need to prove that \( G \cong \mathbb{Z}_4 \) does not admit an O2SR of valency 3. From the discussion in Case 1.2 of Lemma 3.2, we know that \( T_{0,0}, T_{1,1} \subseteq \{\{x\}, \{x^3\}\} \) and \( T_{1,0}, T_{0,1} \subseteq \{\{1, x\}, \{x, x^2\}\} \) with \( T_{1,0} \neq T_{0,1} \). 

Without loss of generality, assume \( T_{0,0} = T_{1,1} = \{x\} \), \( T_{0,1} = \{1, x\} \), and \( T_{1,0} = \{x, x^2\} \). Then we obtain \( \Gamma = \mathrm{Cay}(G_i, T_{i,j} : i, j \in \mathbb{Z}_2) \), which consists of the following arcs:
\[
\begin{aligned}
&1_0 \to x_0,\ 1_0 \to 1_1,\ 1_0 \to x_1,\ x_0 \to x^2_0,\ x_0 \to x_1,\ x_0 \to x^2_1, \\
&x^2_0 \to x^3_0,\ x^2_0 \to x^2_1,\ x^2_0 \to x^3_1,\ x^3_0 \to 1_0,\ x^3_0 \to x^3_1,\ x^3_0 \to 1_1, \\
&1_1 \to x_1,\ 1_1 \to x_0,\ 1_1 \to x^2_0,\ x_1 \to x^2_1,\ x_1 \to x^2_0,\ x_1 \to x^3_0, \\
&x^2_1 \to x^3_1,\ x^2_1 \to x^3_0,\ x^2_1 \to 1_0,\ x^3_1 \to 1_1,\ x^3_1 \to 1_0,\ x^3_1 \to x_0.
\end{aligned}
\]

Substituting \( \Gamma \) into the aforementioned program, we obtain all automorphisms of \( \Gamma \)  as shown in the following figure. Clearly, from the Mathematica results, we have \( |\mathrm{Aut}(\Gamma)| = 8 > |G| = 4 \), so it follows that \( \mathrm{Aut}(\Gamma) \not\cong G \).

\begin{figure}[H]
  \centering
  \includegraphics[width=1.1\linewidth]{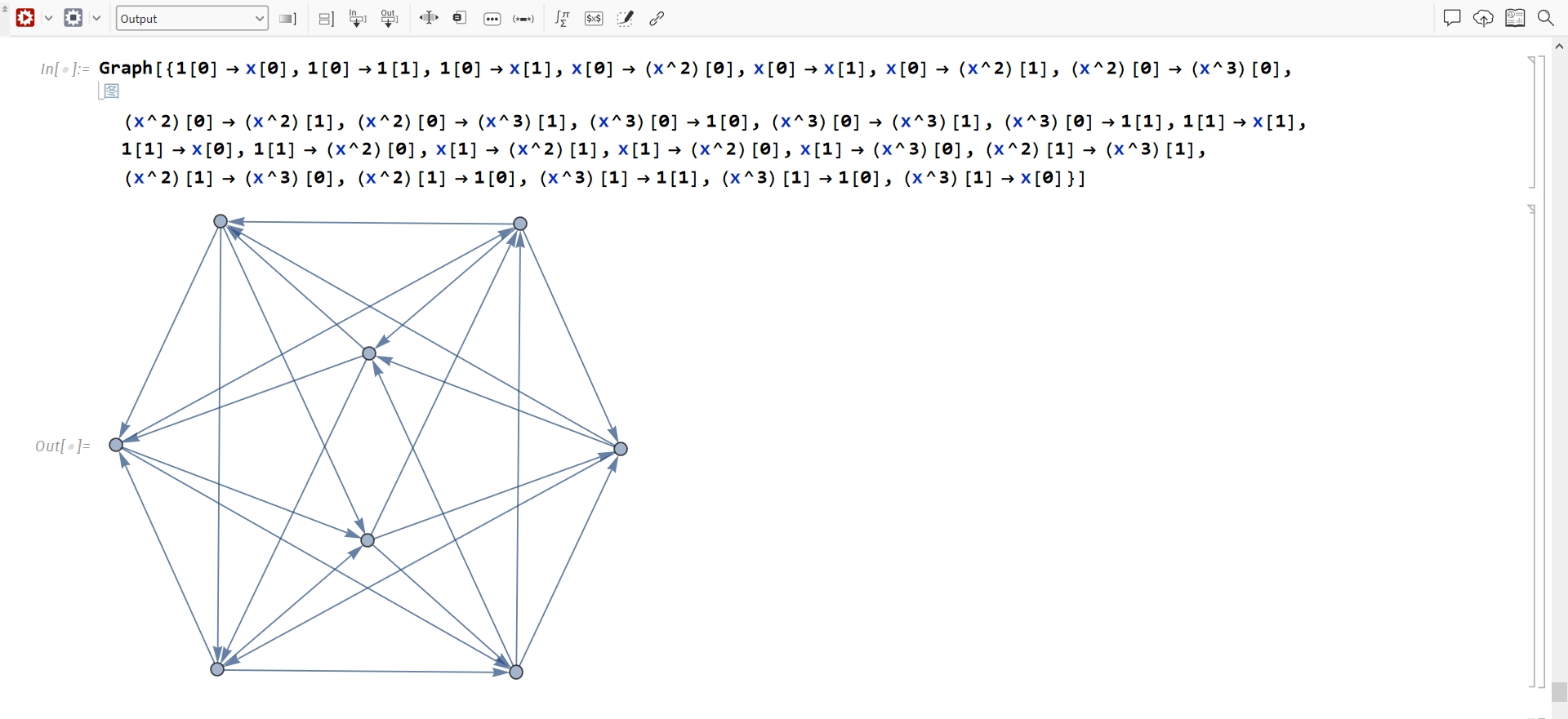}
\end{figure}

\begin{figure}[H]
  \centering
  \includegraphics[width=1.1\linewidth]{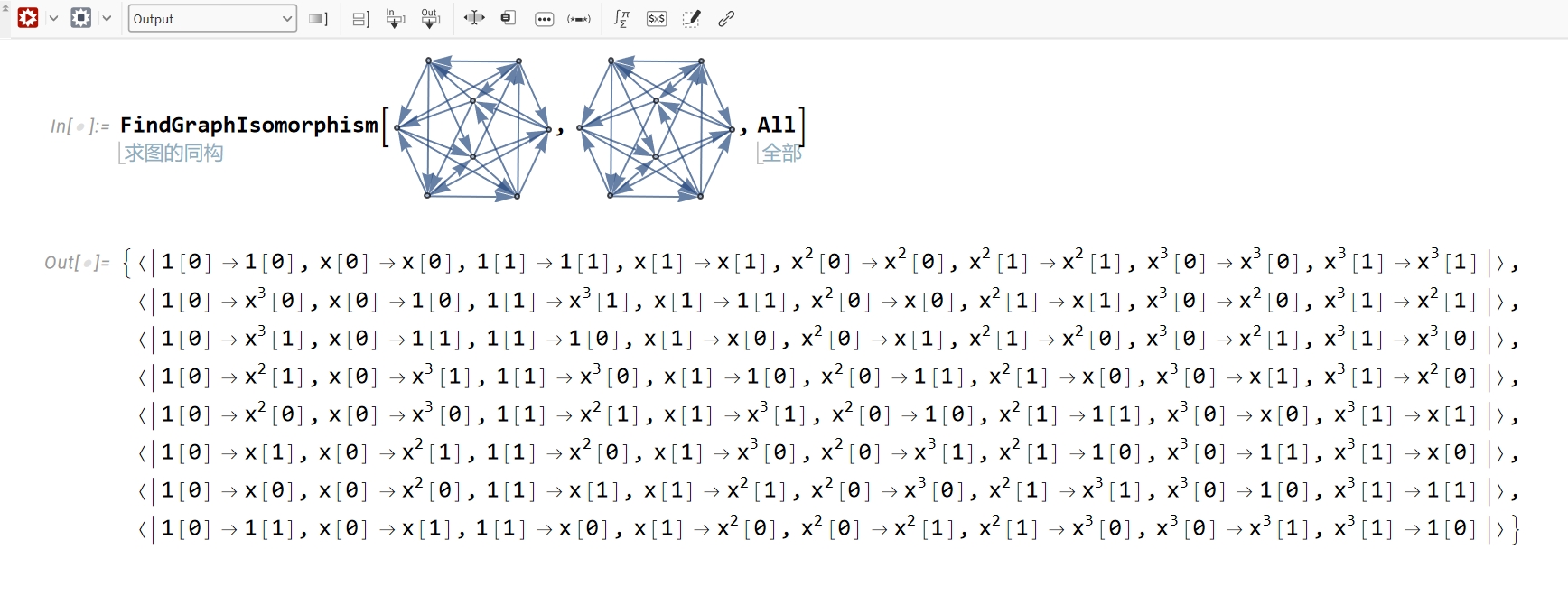}
\end{figure}

\section{Acknowledgments}

We gratefully acknowledge the support of the Graduate Innovation Program of China University of Mining and Technology (2025WLKXJ146), the Fundamental Research Funds for the Central Universities, and the Postgraduate Research and Practice Innovation Program of Jiangsu Province (KYCX25\_2858) for this work.

\end{document}